\documentclass[11pt]{article}

\usepackage{color}
\usepackage{latexsym}
\usepackage{dsfont}
\usepackage{amssymb}
\usepackage{comment}
\usepackage{graphicx}
\usepackage{amsmath, amsfonts,amssymb,theorem,euscript,array,enumerate,amsfonts,mathrsfs}
\usepackage{hyperref}
\usepackage{appendix}

\usepackage{stmaryrd}

\usepackage{mathtools}
\mathtoolsset{showonlyrefs}

\usepackage{comment}

\usepackage{tikz}
\usetikzlibrary{fit,matrix,chains,positioning,decorations.pathreplacing,arrows}

\usepackage{mathtools}
\mathtoolsset{showonlyrefs}

\usepackage[a4paper,margin=1in,heightrounded]{geometry}

\usepackage{algpseudocode}
\usepackage{algorithm}

\usepackage[normalem]{ulem}

\newcommand{\ba}{{\bf a}}

\def \b1{\bf{1}}

\def \R{\mathbb{R}}

\def \E{\mathbb{E}}

\def \P{\mathbb{P}}

\def\esssup_#1{\underset{#1}{\mathrm{ess\,sup\, }}}

\def\argmin_#1{\underset{#1}{\mathrm{argmin\, }}}
\def\argmax_#1{\underset{#1}{\mathrm{argmax\, }}}

\def \Gc{{\cal G}}

\def \Rc{{\cal R}}

\def \Nc{{\cal N}}

\def \Vc{{\cal V}}
\def \Wc{{\cal W}}

\def \ep{\hbox{ }\hfill$\Box$}

\def\bx{{\boldsymbol x}}

\def\bm{{\bf m}}

\def\beqs{\begin{eqnarray*}}
\def\enqs{\end{eqnarray*}}
\def\beq{\begin{eqnarray}}
\def\enq{\end{eqnarray}}

\newtheorem{Theorem}{Theorem}[section]
\newtheorem{Definition}{Definition}[section]
\newtheorem{Proposition}{Proposition}[section]

\newtheorem{Lemma}{Lemma}[section]

\numberwithin{equation}{section}

\begin{document}


\author{
M\'ed\'eric MOTTE
}

\title{Gaussian Cumulative Prospect Theory}

\maketitle

\abstract{
We propose a novel parametrization of Cumulative Prospect Theory (CPT), as developed by Daniel Kahneman and Amos Tversky, that yields an explicit gamble valuation formula for Gaussian reward distributions. Specifically, we define parametric functions \( v_{\theta} \), \( w^{-}_{\theta} \), and \( w^{+}_{\theta} \) satisfying three key properties. The first, \emph{validity}, ensures that for any parameter \(\theta\), the functions conform to the qualitative principles of CPT: \( v_{\theta} \) is concave over gains and convex over losses with a steeper slope for losses; \( w^{-}_{\theta} \) and \( w^{+}_{\theta} \) are increasing, exhibit inverse S-shaped curves, and map 0 to 0 and 1 to 1. The second, \emph{richness}, guarantees that the parametrization is expressive enough to capture a wide range of behaviors: \( v_{\theta} \) can exhibit arbitrary asymptotic behavior and convergence rates, while \( w^{-}_{\theta} \) and \( w^{+}_{\theta} \) can achieve any specified crossover points and slopes. The third, \emph{explicit valuation}, ensures that for any \(\theta\), the CPT valuation of a Gaussian-distributed gamble (with arbitrary mean and variance) can be computed in closed form---enabling efficient approximations for bell-shaped reward distributions. This framework is designed for scalable and rapid computation, making it particularly suited for applications involving large populations. We demonstrate its practicality through two illustrative examples in population-level CPT modeling.


}


\section{Introduction}
In this chapter, we introduce a parametric model for Cumulative Prospect Theory (CPT), originally developed by Daniel Kahneman and Amos Tversky~\cite{Kahneman:1979aa}. Our formulation also implicitly encompasses one of CPT's foundational inspirations, namely John Quiggin's rank-dependent expected utility theory~\cite{Quiggin:1982aa}, although we do not address the latter explicitly in this work.

The central advantage of our model lies in its ability to yield an analytical expression for the valuation of gambles with Gaussian reward distributions. This is particularly appealing given that many real-world random variables approximate Gaussian behavior, as they often arise from the aggregation of numerous small, independent effects—a phenomenon well explained by the Central Limit Theorem. Consequently, our approach is applicable to a broad class of practical decision-making scenarios involving bell-shaped reward distributions.

Having an explicit valuation formula offers several key benefits. First, it significantly reduces the computational burden associated with evaluating risky choices, especially in large-scale population studies where non-analytical models would be computationally prohibitive. Second, the availability of a closed-form expression enables the derivation of analytical gradients and Hessians with respect to model parameters. This facilitates the use of efficient gradient-based optimization techniques—such as (stochastic) gradient descent or Newton's method—in applications like marketing or behavioral modeling, where one seeks to optimize outcomes under risk. Importantly, our approach avoids the numerical instabilities often encountered in models that rely on approximate or simulation-based valuation.

\subsection{Origins and motivations}
The two non-expected utility theories that we focus on in this work—Kahneman and Tversky's Cumulative Prospect Theory (CPT, \cite{Kahneman:1979aa}) and John Quiggin's Rank-Dependent Expected Utility (RDEU) theory (\cite{Quiggin:1982aa})—are still widely considered the most compelling alternatives to the classical Expected Utility Theory (EUT) for modeling decision-making under risk.

Both theories serve as alternatives to EUT, though they can also be interpreted as generalizations of it. EUT, originally proposed by Daniel Bernoulli, aims to model the maximum rational price $V(R)$ one should be willing to pay to engage in a gamble $R$, where $R$ is a real-valued random variable. At the time, the intuitive assumption was that $V(R)$ should equal the expected monetary reward, i.e., $V(R) = \mathbb{E}[R]$. However, Bernoulli challenged this idea with the famous St. Petersburg paradox, thereby demonstrating that this formulation could not adequately describe actual human behavior.

To resolve this paradox, Bernoulli introduced the notion of a \emph{utility function}, suggesting that individuals assign subjective value to monetary outcomes. Under this view, the price $V(R)$ someone is willing to pay corresponds to the expectation of the utility of the reward, i.e., $V(R) = \mathbb{E}[v(R)]$. This marked the first formalization of EUT.

EUT gained significant traction in the mid-20th century, particularly following John von Neumann's foundational work in \emph{Games and Economic Behavior} (\cite{Neumann:2007aa}), where he demonstrated that EUT follows from a set of compelling \emph{axioms of rationality}. Despite its theoretical elegance, numerous empirical studies (\cite{Lichtenstein:1971aa, Lindman:1971aa}) later revealed systematic deviations from EUT in observed human behavior. Today, a substantial body of evidence confirms that EUT often fails to capture actual choice patterns.

To address these empirical inconsistencies, two major research branches emerged: the \emph{conventional} and the \emph{non-conventional} approaches to decision-making under risk. Notably, some scholars have proposed hybrid models bridging both perspectives (\cite{Rubinstein:1988aa, Darity:1992aa, Tversky:1992aa}).

The non-conventional path led to the development of Prospect Theory by Kahneman and Tversky—a theory that earned a Nobel Prize and emphasized psychological realism over axiomatic rigor (\cite{Kahneman:1979aa}). Meanwhile, the conventional approach culminated in Quiggin’s RDEU theory (\cite{Quiggin:1982aa}), which was paradoxically inspired by the non-conventional Prospect Theory itself.

RDEU attracted substantial academic interest (\cite{Machina:1982aa}). It has been thoroughly axiomatized (\cite{Segal:1990aa, Wakker:1994aa, Abdellaoui:2000aa, Yaari:1987aa, Wakker:1994ab}), extended and generalized (\cite{Chew:1989aa, Green:1988aa}), and further explored in various theoretical and empirical directions (\cite{Starmer:1989aa, Luce:1991aa, Tversky:1992aa}).

Subsequently, Kahneman and Tversky refined their original theory into Cumulative Prospect Theory (CPT, \cite{Kahneman:1979aa}), which combined the core insights of Prospect Theory with formal elements from RDEU theory, including the use of probability weighting functions applied to cumulative probabilities.

\vspace{1mm}

A theory of choice under risk—such as Bernoulli’s EUT, Quiggin’s RDEU, or Kahneman and Tversky’s CPT—models how, given a finite set of \emph{prospects} (i.e., a collection of real-valued random variables or gambles $\Gc$), an individual assigns a deterministic value $V(R)$ to each gamble $R \in \Gc$, and then selects the gamble with the highest valuation:
\[
R_\star = \arg\max_{R \in \Gc} V(R).
\]
Such models are highly relevant for predicting individual decisions, and thus have important applications in marketing, finance, policy-making, and behavioral economics.

Choice theories find natural applications in domains such as advertising (predicting whether an individual will click on an advertisement and purchase a product) and political campaigns (forecasting voting behavior). These applications typically require modeling decisions across large populations—web users, consumers, or citizens—which in turn demands the evaluation of gamble valuations for millions or billions of individuals. Since gamble valuations under classical models are defined via integrals, relying on numerical methods for each evaluation can become computationally prohibitive and imprecise at scale. To address this, the goal of this paper is to propose a version of Cumulative Prospect Theory (CPT) specifically designed to yield explicit analytical formulas for gamble valuations, thereby enabling faster, more accurate, and scalable decision modeling.

\subsection{Cumulative Prospect Theory}\label{cpt}

Cumulative Prospect Theory (CPT) assigns a deterministic value to a random prospect based on psychologically grounded transformations of outcomes and probabilities. The evaluation of a prospect $R$ by an individual is defined through two essential components:

\begin{enumerate}
    \item \textbf{Functions with constraints:} There exist three functions:
    \begin{itemize}
        \item a \emph{value function} $v: \R \to \R$,
        \item two \emph{weighting functions} $w^- , w^+ : [0,1] \to [0,1]$,
    \end{itemize}
    satisfying the following empirical constraints, denoted as \textbf{(C)}:
    \begin{itemize}
        \item $v(0)=0$; $v$ is concave on $\R_+$, convex on $\R_-$, with steeper curvature for losses than for gains,
        \item $w^-(0)=w^+(0)=0$, $w^-(1)=w^+(1)=1$; both $w^-$ and $w^+$ are increasing inverse-$S$ shaped functions.
    \end{itemize}

    \item \textbf{Gamble valuation formula:} The subjective value attributed to the prospect $R$ is given by:
    \begin{align}
        \label{eq-valuation}
        V(R) = \int_{-\infty}^0 v(r) \, d(w^- \circ F_R)(r) + \int_{0}^{+\infty} v(r) \, d(w^+ \circ \bar{F}_R)(r),
    \end{align}
    where $F_R$ is the cumulative distribution function (CDF) of $R$, and $\bar{F}_R$ is its tail distribution function.
\end{enumerate}

Alternatively, the gamble valuation $V(R)$ can be described using a probabilistic perspective as follows:
\begin{enumerate}
    \item The random variable $R$ is decomposed into gains and losses: $R^+ := \max(R,0)$, $R^- := \min(R,0)$, so that $R = R^+ + R^-$ almost surely.
    \item The cumulative distribution of $R^-$ is distorted by applying $w^-$, yielding a new distorted distribution. Let $\tilde{R}^-$ be a random variable with this distorted CDF.
    \item Similarly, the tail distribution of $R^+$ is distorted by applying $w^+$, yielding a distorted tail function. Let $\tilde{R}^+$ be a random variable with this distorted tail function.
    \item The valuation is then computed as:
    \[
    V(R) = \mathbb{E}[v(\tilde{R}^-)] + \mathbb{E}[v(\tilde{R}^+)].
    \]
\end{enumerate}

\bigskip

It is important to emphasize that both components of CPT—the \textbf{constraints (C)} and the \textbf{valuation formula}—are equally fundamental. The constraints on $v$, $w^-$, and $w^+$ are not arbitrary: they result from extensive experimental evidence gathered by Kahneman and Tversky. Thus, any modeling or application of CPT must faithfully incorporate these empirical properties.

In addition to the general shape constraints, several analytical parametrizations of the weighting functions $w^\pm$ have been proposed to fit empirical data. Examples include:

\begin{enumerate}
    \item \cite{Tversky:1992aa} proposes a parametrization based on a parameter $\gamma \in [0,1]$:
    \[
    w_\gamma(p) = \frac{p^\gamma}{(p^\gamma + (1-p)^\gamma)^{1/\gamma}}, \quad \forall p \in [0,1].
    \]
    
    \item \cite{Prelec:1998ab} introduces the function:
    \[
    w_\gamma(p) = \exp\left(-(-\ln p)^\gamma\right), \quad \forall p \in [0,1].
    \]

    \item \cite{Tversky:1995aa} suggests the log-odds distortion function $w_{p_0,\gamma}$, defined via the identity:
    \begin{align}
        \label{eq-logodds}
        \mathrm{Lo}(w_{p_0,\gamma}(p)) = \gamma \, \mathrm{Lo}(p) + (1-\gamma) \, \mathrm{Lo}(p_0), \quad \forall p \in (0,1),
    \end{align}
    where $\mathrm{Lo}(p) := \ln\left(\frac{p}{1-p}\right)$ denotes the log-odds transformation.
\end{enumerate}

All three functions exhibit the characteristic inverse-$S$ shape, although this is more easily seen graphically than proved analytically.

\bigskip

According to \cite{Luce}, these parametrizations are difficult to distinguish empirically. Consequently, the specific form of the weighting function is often less critical, provided that:
\begin{itemize}
    \item The functions retain an inverse-$S$ shape, and
    \item The class is sufficiently expressive to approximate general inverse-$S$ curves.
\end{itemize}

\bigskip

As CPT was originally rooted in experimental psychology, the above parametrizations have mostly been used for \emph{fitting empirical data}, particularly to model the behavior of small populations in controlled environments. In such contexts, exact analytical formulas are often unnecessary; discretization and numerical techniques suffice to obtain accurate empirical fits.

\bigskip

However, CPT also has potential applications beyond experimental psychology, such as in commercial or political problems involving large populations. In these domains, the goal is not merely to fit data but to \emph{predict behavior} (e.g., estimating the proportion of individuals making a certain choice) or to \emph{optimize outcomes} (e.g., designing a product or policy that maximizes adoption or votes). Since these applications involve evaluating large numbers of prospects, repeated computation of integrals like those in \eqref{eq-valuation} can become prohibitively expensive.

In such cases, developing analytical approximations or closed-form solutions for the CPT valuation could significantly reduce computational costs and enable real-time decision-making at scale. 

\subsection{Contributions of this work}
Our main contribution is to propose a parametrization of Cumulative Prospect Theory (CPT). Specifically, we introduce parametrized families of:

\begin{enumerate}
    \item reward distributions $\Rc$, from which the random variable $R$ (the gamble) is drawn,
    \item value functions $\Vc$, from which the function $v$ is drawn,
    \item weighting functions $\Wc$, from which the distortion functions $w^-$ and $w^+$ are drawn,
\end{enumerate}

such that each $v \in \Vc$ and $w \in \Wc$ satisfies the theoretical constraints discussed above (point 1 in Section \ref{cpt}), while being flexible enough to approximate any function obeying these constraints. Crucially, this parametrization yields a closed-form valuation formula, making it well-suited for scenarios requiring the valuation of a large number of gambles—such as in large-scale population problems. The proposed classes are as follows:

\begin{enumerate}
    \item[1.] \textbf{Reward Distributions $\Rc$:} The family of Gaussian distributions, parametrized by their mean $\mu$ and standard deviation $\sigma$.

    \item[2.] \textbf{Value Functions $\Vc$:} The class of piecewise exponential utility functions given by:
    \begin{align*}
    v_{m^-, V^-, a^-, m^+, V^+, a^+}(x) &= 
    \begin{cases}
        -\left(m^- x + V^-(1 - e^{-a^-(-x)})\right), & \text{if } x < 0, \\
        \phantom{-}\left(m^+ x + V^+(1 - e^{-a^+ x})\right), & \text{if } x \geq 0,
    \end{cases}
    \end{align*}
    with parameter constraints $m^- \geq m^+$, $V^- \geq V^+$, and $a^- \geq a^+$ to reflect empirical CPT findings.

    \item[3.] \textbf{Weighting Functions $\Wc$:} The class defined by
    \begin{align*}
    w_{p_0, \gamma}(p) = \Nc\left(\gamma \Nc^{-1}(p) + (1 - \gamma) \Nc^{-1}(p_0)\right), \quad \forall p \in [0,1],
    \end{align*}
    where $\Nc$ denotes the cumulative distribution function of the standard normal distribution. This is analogue to the log-odds distortion function \eqref{eq-logodds}, but replacing $\mathrm{Lo}^{-1}$ by $\Nc^{-1}$. We show that this family still yields inverse-$S$ shaped weighting functions (analytically provable), while enabling tractable computations.
\end{enumerate}

Our main theoretical result is the following:

\begin{Theorem}
Let $v \in \Vc$ and $w^- = w_{p_0^-, \gamma^-},\ w^+ = w_{p_0^+, \gamma^+} \in \Wc$. Then:
\begin{itemize}
    \item \textbf{Validity:} The functions $v$ and $w^\pm$ satisfy the CPT constraints \textbf{(C)}.
    
    \item \textbf{Richness:} The families $\Rc$, $\Vc$, and $\Wc$ are expressive: $\Rc$ includes all Gaussian distributions; $\Vc$ allows control over asymptotes and rates of convergence to these asymptotes; $\Wc$ supports arbitrary crossover points and slopes at the crossover.

    \item \textbf{Analytic Valuation Formula:} For a Gaussian gamble with mean $\mu$ and standard deviation $\sigma$, the CPT valuation is given explicitly by:
    \begin{align*}
    &V_{\mu,\sigma, p^-_0,\gamma^-, p^+_0,\gamma^+, m^-,V^-,a^-,m^+,V^+,a^+} \\
    &= -\bm^- \left( \bar{\bx} \Nc(\bar{\bx}) - \frac{1}{\sqrt{2\pi}} e^{- \frac{1}{2} \bar{\bx}^2 } \right) 
    - V^- \left( \Nc(\bar{\bx}) - e^{- \ba^- \bar{\bx} + \frac{(\ba^-)^2}{2}} \Nc(\bar{\bx} - \ba^-) \right) \\
    &\quad + \bm^+ \left( \bx \Nc(\bx) - \frac{1}{\sqrt{2\pi}} e^{- \frac{1}{2} \bx^2 } \right) 
    + V^+ \left( \Nc(\bx) - e^{- \ba^+ \bx + \frac{(\ba^+)^2}{2}} \Nc(\bx - \ba^+) \right),
    \end{align*}
    where
    \begin{align*}
    \hat{\mu} &= \mu - \sigma (\gamma^{-1} - 1) \Nc^{-1}(p_0), \\
    \bar{\hat{\mu}} &= \mu + \sigma (\gamma^{-1} - 1) \Nc^{-1}(p_0), \\
    \hat{\sigma} &= \sigma \gamma^{-1}, \\
    \bx &= \frac{\hat{\mu}}{\hat{\sigma}}, \quad \bar{\bx} = \frac{\bar{\hat{\mu}}}{\hat{\sigma}}, \\
    \bm^\pm &= \hat{\sigma} m^\pm, \quad \ba^\pm = \hat{\sigma} a^\pm.
    \end{align*}
\end{itemize}
\end{Theorem}

We then illustrate the practical relevance of the proposed parametrization by discussing natural examples of large population problems, where the need to compute a vast number of gamble valuations arises and where the analytical formulas derived in this paper lead to drastic reductions in computational cost and time.

\section{The Model}
In this section, we introduce and motivate each parametric class of functions used in our parametric CPT model.

\subsection{The Class $\Rc$ of Reward Probability Distributions}
A probability distribution that is both ubiquitous in Nature and mathematically tractable is the Gaussian (or Normal) distribution. Therefore, we assume that the class of rewards $\Rc$ consists of Gaussian random variables with arbitrary mean and variance, i.e.,
\begin{equation}
\Rc := \left\{ \Nc(\mu, \sigma),\ \mu \in \R,\ \sigma \in \R_+ \right\}.
\end{equation}

More generally, one could consider distributions that are only approximately Gaussian with the same mean and variance. This is frequently observed in Nature and is justified mathematically by the Central Limit Theorem. However, to maintain simplicity and tractability, we restrict ourselves to the case where rewards $R$ follow a Gaussian distribution exactly.

\subsection{The Class $\Wc$ of Weighting Functions}
We now define the class of weighting functions $w$ used in the model. These functions must satisfy two main criteria: (1) they must exhibit an increasing inverse $S$-shaped curve, and (2) they should interact harmoniously with Gaussian-distributed rewards, as assumed above.

\vspace{3mm}

We begin by formalizing the necessary conditions for a function $w : [0,1] \to [0,1]$ to be a valid weighting function under CPT.

\begin{Definition}[Valid Weighting Function for CPT]
A valid weighting function is a differentiable function $\phi: [0,1] \to [0,1]$ satisfying:
\begin{itemize}
    \item $\phi(0) = 0$, $\phi(1) = 1$,
    \item $\phi'(p) > 0$ for all $p \in [0,1]$ (monotonicity),
    \item $\phi$ has a unique inflection point $p^\star \in (0,1)$ such that $\phi''(p^\star) = 0$, $\phi''(p) < 0$ for $p \in [0, p^\star)$, and $\phi''(p) > 0$ for $p \in (p^\star, 1]$ (inverse $S$-shape).
\end{itemize}
\end{Definition}

We now introduce a parametric family of probability distortion functions tailored to our model. Inspired by the log-odds weighting function used in \cite{Tversky:1995aa}, defined by:
\begin{equation}
\mathrm{Lo}(w_{p_0,\gamma}(p)) = \gamma \mathrm{Lo}(p) + (1-\gamma) \mathrm{Lo}(p_0), \quad \forall p \in [0,1],
\end{equation}
we observe that it is based on the log-odds (quantile) function of the logistic distribution, whose shape resembles that of the Gaussian distribution.

However, since we assume normally distributed rewards, it is more natural and convenient to replace the logistic quantile with the Gaussian quantile. This leads to the following definition:

\begin{Definition}[Normal Distortion Function]
The \emph{normal distortion function} $w^N_{p_0, \gamma}: [0,1] \to [0,1]$ with parameters $p_0 \in [0,1]$, $\gamma \in [0,1]$ is defined as:
\begin{equation}
w^N_{p_0, \gamma}(p) := \Nc\left( \gamma \Nc^{-1}(p) + (1 - \gamma) \Nc^{-1}(p_0) \right), \quad \forall p \in [0,1],
\end{equation}
where $\Nc$ denotes the standard normal CDF and $\Nc^{-1}$ its quantile function.
\end{Definition}

\subsubsection{Validity of $w^N_{p_0, \gamma}$ for $p_0 \in [0,1]$, $\gamma \in [0,1]$}
We now show that $w^N_{p_0, \gamma}$ satisfies the criteria of a valid weighting function.

\begin{Proposition}\label{prop-normal}
For all $p_0, \gamma \in [0,1]$, the normal distortion function $w^N_{p_0, \gamma}$ is a valid weighting function for CPT.
\end{Proposition}

{\bf Proof.}
Let us define $u_0 := (1 - \gamma)\Nc^{-1}(p_0)$. Then,
\[
w^N_{p_0, \gamma}(p) = \Nc\left( \gamma \Nc^{-1}(p) + u_0 \right), \quad \forall p \in [0,1].
\]
We first verify the boundary conditions:
\[
w^N_{p_0, \gamma}(0) = \Nc(-\infty) = 0, \quad w^N_{p_0, \gamma}(1) = \Nc(+\infty) = 1.
\]

Next, compute the first derivative:
\[
(w^N_{p_0, \gamma})'(p) = \gamma (\Nc^{-1})'(p) \cdot n\left( \gamma \Nc^{-1}(p) + u_0 \right),
\]
\beqs 
w'_{p_0,\gamma}(p)=\gamma (N^{-1})'(p)N'(\gamma N^{-1}(p) + u_0)=\gamma \frac{n(\gamma N^{-1}(p) + u_0)}{n(N^{-1}(p))}>0,
\enqs 
where $n(x) = \frac{1}{\sqrt{2\pi}}e^{-x^2/2}$ is the standard normal density, and where the inequality simply comes from the positivity of the gaussian density. Let us now establish the existence of a unique inflection point on $[0,1]$. Define $f:[0,1]\rightarrow \R$ by $f(p)=\ln(w'_{p_0,\gamma}(p))$, for all $p\in[0,1]$. Note that $f'=\frac{w''_{p_0,\gamma}}{w'_{p_0,\gamma}}$, and thus, $w''_{p_0,\gamma}(p)$ has only one zero if and only if $f'$ has only one zero. We have
\beqs 
f(p)=\ln(\gamma) +\frac{1}{2}\Big( N^{-1}(p)^2-(\gamma N^{-1}(p)+u_0)^2\Big)
\enqs 
Then,
\beqs 
f'(p)&=&\frac{1}{2}\Big( 2(N^{-1})'(p)N^{-1}(p)-2\gamma(N^{-1})'(p)(\gamma N^{-1}(p)+u_0)\Big)\\
&=&(N^{-1})'(p)\Big( N^{-1}(p)-\gamma(\gamma N^{-1}(p)+u_0)\Big)\\
\enqs 
Note that $(N^{-1})'(p)=\frac{1}{n(N^{-1}(p))}>0$ for all $p\in[0,1]$ by the positivity of $n$. Therefore, we only need to study the sign of $N^{-1}(p)-\gamma(\gamma N^{-1}(p)+u_0)$. We have:
\beqs 
&&N^{-1}(p)-\gamma(\gamma N^{-1}(p)+u_0)\geq 0\\
&&\Leftrightarrow N^{-1}(p)(1-\gamma^2)-\gamma u_0\geq 0\\
&&\Leftrightarrow N^{-1}(p)\geq \frac{\gamma u_0}{1-\gamma^2}\Leftrightarrow p\geq N\Big(\frac{\gamma N^{-1}(p_0)}{1+\gamma}\Big)
\enqs 
which means that $w''_{p_0,\gamma}$ is negative on $[0,N\Big(\frac{\gamma N^{-1}(p_0)}{1+\gamma}\Big)]$ and positive on $[N\Big(\frac{\gamma N^{-1}(p_0)}{1+\gamma}\Big), 1]$.
\ep

\subsubsection{Gaussian stability property}
The class of weighting functions $\mathcal{W}$ possesses a remarkable feature: it preserves the Gaussian structure under distortion. More precisely, applying a weighting functions in $\mathcal{W}$ to the cumulative (or tail) distribution function of a Gaussian random variable yields the cumulative (or tail) distribution function of another Gaussian random variable.

The class of weighting functions $\Wc$ has another particularity: it stabilizes gaussian distributions. More precisely, we have the following result.



\begin{Lemma}\label{lem-gaus}
Let $\Nc_{\mu,\sigma}$ be the distribution function of a gaussian variable with mean $\mu$ and standard deviation $\sigma$, and let $w^N_{p_0,\lambda}\in \mathcal{W}$ be a weighting function with crossover point $p_0$ and slope at the crossover point $\lambda$. Then, the function $w^N_{p_0,\lambda}\circ \Nc_{\mu,\sigma}$ is the distribution function of a gaussian variable with mean $\hat{\mu}:=\mu-\sigma(\gamma^{-1}-1) \Nc^{-1}(p_0)$ and standard deviation $\hat{\sigma}:=\sigma \gamma^{-1}$, and $w^N_{p_0,\lambda}\circ \bar{\Nc}_{\mu,\sigma}$ is the tail function of a gaussian variable with mean $\bar{\hat{\mu}}:=\mu+ \sigma(\gamma^{-1}-1) \Nc^{-1}(p_0)$ and standard deviation $\bar{\hat{\sigma}}:=\sigma \gamma^{-1}=\hat{\sigma}$. 
\end{Lemma}
{\bf Proof.}
We simply have
\beqs 
&&w^N_{p_0,\lambda}\circ \Nc_{\mu,\sigma}(x)= \Nc(\gamma \Nc^{-1}(\Nc(\frac{x-\mu}{\sigma})) + (1-\gamma) \Nc^{-1}(p_0))\\
&=& \Nc\left(\gamma \frac{x-\mu}{\sigma} + (1-\gamma) \Nc^{-1}(p_0)\right)= \Nc\left(\frac{x-( \mu-(\gamma^{-1}-1)\sigma \Nc^{-1}(p_0))}{\sigma \gamma^{-1}} \right)\\
&=&\Nc_{\hat{\mu},\hat{\sigma}}(x),
\enqs 
where 
\beqs 
\hat{\mu}= \mu-\sigma(\gamma^{-1}-1) \Nc^{-1}(p_0), \quad \hat{\sigma}=\sigma \gamma^{-1}
\enqs 
Likewise, 
\beqs 
&&w^N_{p_0,\lambda}\circ \bar{\Nc}_{\mu,\sigma}(x)= w^N_{p_0,\lambda}\circ \Nc_{-\mu,\sigma}(-x)\\
&=&\Nc_{\widehat{-\mu},\hat{\sigma}}(-x)=\bar{\Nc}_{-\widehat{-\mu},\hat{\sigma}}(x)=\bar{\Nc}_{\tilde{\mu},\hat{\sigma}}(x)
\enqs 
where 
\beqs 
\tilde{\mu}=\mu+ \sigma(\gamma^{-1}-1) \Nc^{-1}(p_0)
\enqs 
\ep

Therefore, the class of weighting functions $\Wc$ ``stabilizes'' the set of gaussian distributions.

\subsection{The Value Function}

We now introduce the class $\mathcal{V}$ of value functions in a natural way.

So far, our model has specified two components:
\begin{enumerate}
    \item The class of reward distributions is taken to be the family of Gaussian distributions.
    \item The class of weighting functions is the family of normal weighting functions $\mathcal{W}$.
\end{enumerate}

Thus, when modeling decision-making under risk, a typical instance involves choosing parameters $\mu$ and $\sigma$ for the Gaussian reward distribution associated with a risky choice, and parameters $\lambda$ and $p_0$ for the individual's weighting function. 

By the definition of valuation in Cumulative Prospect Theory (CPT), and invoking Lemma~\ref{lem-gaus}, it follows that the class $\mathcal{V}$ of value functions should ideally be composed of functions that are either explicitly or nearly explicitly integrable with respect to any Gaussian distribution.

What remains is to construct such a value function that also satisfies the essential axioms of Prospect Theory: it must be {\it concave} on $\mathbb{R}_+$, {\it convex} on $\mathbb{R}_-$, {\it steeper} on the negative axis (loss aversion), and {\it vanish at zero}.

\vspace{5mm}

A particularly useful family of functions that integrates well against Gaussian measures is the class of exponential functions of the form $x \mapsto e^{-a x}$. Consider the following value function:
\[
v(x) := 
\begin{cases}
    -m^- x - V^- \left(1 - e^{-a^-(-x)}\right), & \text{if } x < 0, \\
    m^+ x + V^+ \left(1 - e^{-a^+ x}\right), & \text{if } x \geq 0.
\end{cases}
\]
This function is concave on $\mathbb{R}_+$, convex on $\mathbb{R}_-$, and steeper on the negative axis when the parameters satisfy:
\[
V^- > V^+, \quad a^- > a^+, \quad m^- > m^+.
\]

Let us interpret the parameters of this utility function more intuitively. The function is defined piecewise over $\mathbb{R}_-$ and $\mathbb{R}_+$, allowing separate modeling of losses and gains.

Focusing on the positive part ($x \geq 0$), we observe two natural features that allow this function to represent a wide range of increasing concave preferences:
\begin{itemize}
    \item An \textbf{asymptote}, characterizing the behavior as $x \to +\infty$.
    \item A \textbf{rate of convergence} toward this asymptote, describing how quickly the function approaches it.
\end{itemize}

The function $v(x)$ has the asymptotic form:
\[
v(x) \to V^+ + m^+ x \quad \text{as } x \to +\infty,
\]
and for $x \geq 0$, we have the exact expression:
\[
V^+ + m^+ x - v(x) = V^+ e^{-a^+ x}.
\]
This shows that the value function converges exponentially to the asymptote $V^+ + m^+ x$ at rate $a^+$. The function starts at $v(0) = 0$, with an initial vertical offset of $V^+$, and then converges toward its asymptote at an exponential rate governed by $a^+$.

A similar interpretation holds on the negative axis for losses.

\section{The Gamble Valuation Function}

In this section, we derive an explicit expression for the gamble valuation $V(R)$, where the random variable $R \sim \mathcal{N}(\mu, \sigma)$, and where $w^-, w^+ \in \mathcal{W}$ and $v \in \mathcal{V}$.

By ``explicit expression," we mean a formula composed of elementary operations, along with functions that are widely supported and efficiently implemented in standard statistical libraries—namely, the exponential function and the cumulative distribution function (CDF) of the standard normal distribution.

\begin{Theorem}
We have
\beqs 
&&V_{\mu,\sigma, p_0,\lambda, m^-,V^-,a^-,m^+,V^+,a^+}\\
&=&- \bm^- \left(\bx\Nc(\bx)-\frac{1}{\sqrt{2\pi}}e^{-\frac{1}{2}\bx^2}\right)-V^-\left(\Nc(\bx)-e^{-\ba^-\bx+\frac{(\ba^-)^2}{2}}\Nc(\bx-\ba^-)\right)\\
&+&\bm^+ \left(\bx\Nc(\bx)-\frac{1}{\sqrt{2\pi}}e^{-\frac{1}{2}\bx^2}\right)+V^+\left(\Nc(\bx)-e^{-\ba^+\bx+\frac{(\ba^+)^2}{2}}\Nc(\bx-\ba^+)\right)
\enqs 
where, $\bx:=\frac{\hat{\mu}}{\hat{\sigma}}$, $\bar{\bx}:=\frac{\bar{\hat{\mu}}}{\hat{\sigma}}$, $\bm^+:=\hat{\sigma}m^+$, $\bm^-:=\hat{\sigma}m^-$, $\ba^+:=\hat{\sigma}a^+$, $\ba^-:=\hat{\sigma}a^-$, where $\hat{\mu}=\mu-\sigma(\gamma^{-1}-1) \Nc^{-1}(p_0)$, $\bar{\hat{\mu}}=\mu+ \sigma(\gamma^{-1}-1) \Nc^{-1}(p_0)$, $\hat{\sigma}=\sigma \gamma^{-1}$.
\end{Theorem}
\noindent{\bf Proof.}
The idea is to split the computation into several components. First of all, notice that the fact that gains and losses are processed separately (using different parameters for the weighting function and value function), we can clearly focus on computing the part corresponding to the gains: the part a corresponding to losses will take the same form but with different parameters. On $\R_+$, notice that the value function is
\beqs 
v(x)=m^+x+V^+(1-e^{-a^+x}),\quad \forall x\in\R_+.
\enqs 
We can thus separately study the terms $m^+x$, $V^+$, and $-V^+e^{-a^+x}$. The ``gain'' part of gamble $R$, i.e. $R_+$ will have its tail function distorted by $w^+$. By Lemma \ref{lem-gaus}, this will yield another gaussian variable with parameters $\bar{\hat{\mu}}$ and $\hat{\sigma}$ given by Lemma \ref{lem-gaus}. All we are left to do is, given $Z\sim \Nc(\bar{\hat{\mu}},\hat{\sigma})$ to compute separately
\beqs 
\E[m^+Z{\bf 1}_{Z\geq 0}],\quad \E[V^+{\bf 1}_{Z\geq 0}], \quad \E[-V^+e^{-a^+x}{\bf 1}_{Z\geq 0}].
\enqs 
The second term is simply given by $\E[V^+{\bf 1}_{Z\geq 0}]=V^+\Nc(\frac{\bar{\hat{\mu}}}{\hat{\sigma}})$. Let us compute the third term. We have
\beqs 
&&-V^+\E[e^{-a^+Z}{\bf 1}_{Z\geq 0}]=-V^+e^{-a^+\hat{\mu}}\E[e^{-a^+\hat{\sigma}N}{\bf 1}_{N\geq -\frac{\hat{\mu}}{\hat{\sigma}}}]=-V^+e^{-a^+\hat{\mu}}\E[e^{a^+\hat{\sigma}N}{\bf 1}_{N< \frac{\hat{\mu}}{\hat{\sigma}}}]\\
&=&-V^+e^{-a^+\hat{\mu}}\int_{-\infty}^{\frac{\hat{\mu}}{\hat{\sigma}}}e^{a^+\hat{\sigma}x}\frac{1}{\sqrt{2\pi}}e^{-\frac{x^2}{2}}dx=-V^+e^{-a^+\hat{\mu}}e^{\frac{(a^+\hat{\sigma})^2}{2}}\int_{-\infty}^{\frac{\hat{\mu}}{\hat{\sigma}}}\frac{1}{\sqrt{2\pi}}e^{-\frac{(x-a^+\hat{\sigma})^2}{2}}dx\\
&=&-V^+e^{-a^+\hat{\mu}+\frac{(a^+\hat{\sigma})^2}{2}}\P(N+a^+\hat{\sigma} < \frac{\hat{\mu}}{\hat{\sigma}})=-V^+e^{-a^+\hat{\mu}+\frac{(a^+\hat{\sigma})^2}{2}}\P(N < \frac{\hat{\mu}-a^+\hat{\sigma}^2}{\hat{\sigma}})\\
&=&-V^+(e^{-a^+\hat{\mu}+\frac{(a^+\hat{\sigma})^2}{2}}\Nc(\frac{\hat{\mu}-a^+\hat{\sigma}^2}{\hat{\sigma}}))
\enqs 
Finally, the first term is obtained as follows. Notice that
\beqs 
\E[Z{\bf 1}_{Z\geq 0}]&=& -\partial_{a^+=0}\E[e^{-a^+Z}{\bf 1}_{Z\geq 0}]=\hat{\mu}\Nc(\frac{\hat{\mu}}{\hat{\sigma}})-\frac{\hat{\sigma}}{\sqrt{2\pi}}e^{-\frac{1}{2}\left(\frac{\hat{\mu}}{\hat{\sigma}}\right)^2}
\enqs 
\ep

Notice that the explicit expression of $V(R)$ makes it easily derivable in each of the model's parameters. 

 


\section{Applications to Large Population Problems}
In this section, we highlight the practical benefits of the analytical expressions derived for the Gaussian Cumulative Prospect Theory (CPT) model, particularly in the context of large-scale decision-making environments. This work builds upon a line of research investigating prospect-theoretic agents in large populations, including mathematical foundations for modeling such populations \cite{motte2021mathematical}, theoretical mean-field formulations \cite{motte2022mean}, quantitative propagation of chaos results \cite{motte2023quantitative}, and applied models for opinion dynamics \cite{coculescu2024opinion} and marketing optimization \cite{motte2021optimal}. The closed-form expressions derived in the present article enable real-time evaluation and gradient computation in such models, providing a key computational advantage for optimization and equilibrium analysis in large-scale environments.

\vspace{3mm}

Analytical formulas offer numerous advantages: fast computation, interpretability, sensitivity analysis, explicit gradients and Hessians for optimization algorithms (e.g., gradient descent, Newton's method), closed-form solutions to complex problems, real-time plotting, and support for high-frequency decision-making, such as in algorithmic trading. Rather than exhaustively listing all such possibilities, we focus here on two illustrative applications in large population contexts.

\subsection{Optimal Product or Program Design in a Large Population}
While numerical integration is now computationally efficient, it remains significantly slower than evaluating an explicit formula, especially when scaled to millions or billions of evaluations. We present a setting in which repeated computation of gamble valuations naturally arises.

\vspace{2mm}

Consider an agent (e.g., a company or politician) seeking to design a program, such as a marketing or electoral campaign. Let $\mathcal{P}$ denote the set of all possible programs. We consider a population of $N$ individuals, each characterized by a personality trait $e_n \in E$, where $E$ is some personality space. Each individual $e \in E$ has associated model parameters governing their CPT functions, defined by:
\[
p_0^+, \lambda^+, V^+, m^+, a^+, \quad p_0^-, \lambda^-, V^-, m^-, a^-: E \rightarrow \mathbb{R}.
\]
Additionally, for any program $P \in \mathcal{P}$, the outcome distribution perceived by an individual with personality $e$ is Gaussian with mean $\mu(e, P)$ and variance $\sigma(e, P)$. These mappings define a parameter vector $\theta(e, P)$, such that the certainty equivalent $CE_{\theta(e,P)}$ is an explicit function of these parameters. Provided that $\mu$ and $\sigma$ are explicit functions of $P$, the certainty equivalent remains explicit in $P$.

\vspace{2mm}

We assume individuals must make a binary choice (e.g., to purchase---or not purchase---a product, or to vote---or not vote---for a candidate). The decision rule is modeled as:
\[
x_n = \mathbf{1}_{CE_{\theta(e_n, P)} > 0},
\]
i.e. an individual chooses to purchase the product (or vote for the candidate) if they qssign a positive value to this action. Hence, the aggregate proportion of individuals choosing the program $P$ is:
\[
\frac{1}{N} \sum_{n=1}^N \mathbf{1}_{CE_{\theta(e_n, P)} > 0}.
\]
Evaluating this quantity using analytical expressions offers significant computational savings, especially when $N$ is large.

\vspace{2mm}

Furthermore, assume the agent derives a gain $G_N(P)$ from presenting program $P$ to the population, depending on the fraction of individuals selecting it:
\[
G_N(P) = g\left(P, \frac{1}{N} \sum_{n=1}^N \mathbf{1}_{CE_{\theta(e_n, P)} > 0}\right).
\]
The optimization problem becomes:
\[
P^*_N = \arg\max_{P \in \mathcal{P}} G_N(P).
\]
Solving this problem involves repeated evaluation of $G_N(P)$, each requiring $N$ certainty equivalent computations. The difference in computational efficiency between numerical and analytical evaluation becomes dramatic at scale.

\vspace{2mm}

We may also consider the mean-field approximation by replacing the empirical distribution $\frac{1}{N} \sum_{n=1}^N \delta_{e_n}$ with a probability measure $\nu$ on $E$. Letting $\varepsilon \sim \nu$, we define:
\[
G(P) = g\left(P, \mathbb{P}(CE_{\theta(\varepsilon, P)} > 0)\right),
\quad
P^* = \arg\max_{P \in \mathcal{P}} G(P).
\]
Computing $\mathbb{P}(CE_{\theta(\varepsilon, P)} > 0)$ requires sampling multiple realizations of $CE_{\theta(e, P)}$. Again, analytical evaluation offers substantial efficiency. Additionally, gradient-based methods to solve the optimization may rely on computing:
\[
\partial_P \mathbb{P}(CE_{\theta(\varepsilon, P)} > 0)
= \mathbb{E}\left[\partial_P CE_{\theta(\varepsilon, P)} \mid CE_{\theta(\varepsilon, P)} = 0\right],
\]
which further motivates the need for closed-form expressions and derivatives.

\subsection{Equilibrium Computation in a Large-Population Social Game}
We now consider a second application involving social interactions in a large population.

\vspace{2mm}

Assume that a company proposes a fixed product $P$. In contrast to the previous setup, we now posit that the reward experienced by each individual depends not only on the product but also on social dynamics. Specifically, if multiple individuals choose the same option (e.g., buy the product), they derive an additional social reward from the shared experience.

\vspace{2mm}

Let $x = (x_1, \dots, x_N) \in \{0, 1\}^N$ denote the population's choices. Individual $n$'s reward becomes:
\[
R_{P, e_n} + \mathfrak{u}\left(\frac{1}{N} \sum_{i=1}^N x_i\right) 
\sim \mathcal{N}\left(\mu(e_n, P) + \mathfrak{u}\left(\frac{1}{N} \sum_{i=1}^N x_i\right), \sigma(e_n, P)\right),
\]
where $\mathfrak{u}$ is a function modeling the social interaction benefit. The resulting mean now depends on aggregate behavior, and the parameter mapping becomes $\theta(e_n, P, \mathfrak{u}(q))$, where $q = \frac{1}{N} \sum_i x_i$.

\vspace{2mm}

Each individual aims to select $x_n$ such that:
\[
x_n = \mathbf{1}_{CE_{\theta(e_n, P, \mathfrak{u}(q))} > 0}.
\]
As individuals do not know others' choices, this forms a large-population game. Searching equilibriums for this game (eg by iterative algorithms) will then clearly benefit from having closed-form expressions for the certainty equivalent.

\section{Conclusion}
In this paper, we introduced a parametric model for Cumulative Prospect Theory (CPT). Specifically, we defined a parameter space $\Theta$ and a mapping $\theta \in \Theta \mapsto (v_\theta, w^-_\theta, w^+_\theta)$ such that, for any $\theta \in \Theta$, the functions $v_\theta$, $w^-_\theta$, and $w^+_\theta$ constitute \emph{valid} value and probability weighting functions. This construction satisfies two essential properties: a \emph{validity} property, ensuring all resulting functions are admissible within the CPT framework, and a \emph{richness} property, ensuring that any general shape of valid functions can be approximated by a suitable choice of $\theta$.

Additionally, we derived an explicit formula for evaluating Gaussian rewards under CPT preferences, and showed that this valuation is differentiable with respect to any component of the parameter vector $\theta$. We illustrated the practical utility of this analytical expression through two examples involving large populations, where access to closed-form valuations and gradients significantly accelerates computation, often by several orders of magnitude.

These findings suggest promising directions for future research, particularly in the study of large-population decision-making problems. While briefly discussed here, such settings are especially relevant in applications where aggregate behavior matters more than individual actions—for example, in commercial contexts such as product design aimed at maximizing market share, or in political contexts such as shaping policies to secure electoral support.

Our parametric approach thus opens new avenues for applying CPT to real-world problems at scale, bridging theoretical modeling with practical impact.

\bibliographystyle{plain}
\bibliography{bibliography}
 
\end{document}